 \UseRawInputEncoding

\documentclass[12pt]{amsart}
\usepackage{latexsym, amssymb, amsmath, amscd}
 \usepackage{eucal}

    \newtheorem{rema}{Remark}[section]
    \newtheorem{propo}[rema]{Proposition}

   \newtheorem{theo}[rema]{Theorem}
   \newtheorem{def-theo}[rema]{Definition-Theorem}
 
   \newtheorem{defi}[rema]{Definition}
    \newtheorem{lemma}[rema]{Lemma}
    \newtheorem{corol}[rema]{Corollary}
     \newtheorem{exam}[rema]{Example}
  \newtheorem{rmk}[rema]{Remark}

\newtheorem{oprob}[rema]{Open Problem}

	\newcommand{\nno}{\nonumber}

	\newcommand{\p}{\partial}

 \newcommand{\pf}{{\it Proof:}\hspace{2ex}}
 
 \newcommand{\epfv}{\hspace{1em}$\Box$\vspace{1em}}

\newcommand{\bC}{{\mathbb C}}
\newcommand{\bZ}{{\mathbb Z}}
\newcommand{\bR}{{\mathbb R}}

\newcommand{\bN}{{\mathbb N}}

\newcommand{\cZ}{{\mathcal Z}}

\newcommand{\cW}{{\mathcal W}} 

\newcommand{\cA}{{\mathcal A}}

\newcommand{\rad}{{\frak r}}
\newcommand{\nil}{\mbox{\rm nil\,}}
\newcommand{\ann}{\mbox{\rm Ann\,}}

\newcommand{\cB}{{\mathcal B}}

\newcommand{\cE}{{\mathcal E}}

\newcommand{\Ker}{\operatorname{Ker}}
\newcommand{\im}{\operatorname{Im\,}}

\newcommand{\ad}{\operatorname{ad}}

\newcommand{\I}{{\rm I }}

\title[The Radical of the Kernel of a Differential Operator]
{The Radical of the Kernel of a Certain Differential Operator  
and Applications to Locally Algebraic Derivations}

  \author{Wenhua Zhao}      
    \date{August 08, 2022}
\address{Department of Mathematics, Illinois State University, Normal, IL 61761. Email: wzhao@ilstu.edu}

\begin{document}

\begin{abstract}
Let $R$ be a commutative ring, $\mathcal A$ 
an $R$-algebra (not necessarily commutative) and $V$ an $R$-subspace or 
$R$-submodule of $\mathcal A$. By the {\it radical} of $V$  
we mean the set of all elements $a\in \mathcal A$ such that $a^m\in V$ for all $m\gg 0$. We derive (and show) some necessary conditions satisfied by the elements in the radicals of the kernel of some   
(partial) differential operators, such as all differential operators of commutative algebras; the differential operators $P(D)$ of (noncommutative) $\mathcal A$ with certain conditions, 
where $P(\cdot)$ is a polynomial in $n$ commutative free variables and 
$D=(D_1, D_2, \dots, D_n)$ are either commuting locally finite   
$R$-derivations or commuting $R$-derivations of $\mathcal A$ 
such that for each $1\le i\le n$, $\mathcal A$ can be decomposed 
as a direct sum of the generalized eigen-subspaces of $D_i$; etc.    
In particular, we show that the kernel of 
certain differential operators of $\mathcal A$ is a 
Mathieu subspace (see \cite{GIC, MS}) of $\mathcal A$. 
We then apply some results above to study $R$-derivations 
of $\mathcal A$, which are locally algebraic or locally integral 
over $R$.  In particular, we show that if $R$ is an integral domain of characteristic zero and $\mathcal A$ is reduced and torsion-free as an $R$-module, then $\mathcal A$ has no nonzero locally algebraic $R$-derivations. We also show a formula for the determinant of a  
differential vandemonde matrix over a commutative algebra $\mathcal A$. 
This formula not only provides some information for the elements in the radical of the kernel of all ordinary differential operators of $\mathcal A$, but also is interesting on its own right.
\end{abstract}

   









%
%



\thanks{The author has been partially supported 
by the Simons Foundation grant 278638}

 \bibliographystyle{alpha}
    \maketitle


\renewcommand{\theequation}{\thesection.\arabic{equation}}
\renewcommand{\therema}{\thesection.\arabic{rema}}
\setcounter{equation}{0}
\setcounter{rema}{0}
\setcounter{section}{0}

\section{\bf Background and Motivation}\label{S1}
Let $R$ be a commutative ring and $\cA$ 
an $R$-algebra (not necessarily commutative). 
A {\it derivation} $D$ of $\cA$ is  
a map from $\cA$ to $\cA$ such that $D(a+b)=D(a)+D(b)$ and  
$D(ab)=D(a)b+aD(b)$ 
for all $a, b\in \cA$.
If $D$ is also $R$-linear, 
we call it an {\it $R$-derivation}  
of $\cA$. 

For each $a\in \cA$, denote by $\ell_a$ 
the map from $\cA$ to $\cA$ that 
maps $b\in \cA$ to $ab$. We call the associative 
algebra generated by $\ell_a$ ($a\in \cA$) and 
all derivations of $\cA$ {\it the Weyl algebra} 
of $\cA$, and denote it by $\cW(\cA)$. 
The subalgebra of $\cW(\cA)$ generated 
by $\ell_a$ ($a\in R$) and all $R$-derivations 
of $\cA$ will be denoted by $\cW_R(\cA)$. We call 
elements of $\cW(\cA)$ the {\it differential operators} 
of $\cA$. 

For each $\Phi\in\cW(\cA)$, it is well-known and also easy to check 
that there exist some derivations $D=(D_1, D_2, \dots, D_n)$ of $\cA$ 
and a polynomial $P(\xi)\in\cA\langle \xi\rangle$ (the polynomial algebra over $\cA$ in $n$ noncommutative free variables 
$\xi=(\xi_1, \xi_2, \dots, \xi_n)$) such that $\Phi=P(D)$, where 
$P(D)$ throughout this paper is defined by first {\it writing all the coefficients of 
$P(\xi)$ on the most left of the monomials in $\xi$}, and then replacing $\xi_i$ by $D_i$ for all 
$1\le i\le n$. Furthermore, if $\Phi\in \cW_R(\cA)$, the same is true 
with $D_i$ $(1\le i\le n)$ being $R$-derivations of $\cA$ 
and $P(\xi)\in R\langle \xi\rangle$. We call the differential operator 
$\Phi=P(D)$ an {\it ordinary differential operator} of $\cA$,  
if $P(\xi)$ is univariate, and a 
{\it partial differential operator} of $\cA$ if  
$P(\xi)$ is multivariate. 

%

Next, we recall the following two notions of associative algebras that were first introduced in \cite{GIC, MS}.  

\begin{defi} \label{Def-MS}
An $R$-subspace (or $R$-submodule) $V$ of an $R$-algebra $\cA$ is said to  be a {\it Mathieu subspace} (MS)  
of $\cA$ if for all $a, b, c\in \cA$ with 
$a^m\in V$ for all $m\ge 1$, we have 
$ba^mc \in V$ for all $m\gg 0$, i.e., there exists $N\in\bN$ (depending on $a, b, c$) such that $ba^mc \in V$ for all $m\ge N$.
\end{defi}
 
Note that a MS is also called a {\it Mathieu-Zhao space} 
in the literature (e.g., see \cite{DEZ, EN, EH, EKC}, etc.), 
as first suggested by A. van den Essen \cite{E2}. 
The introduction of this notion  
is mainly motivated by the studies in \cite{Ma, IC} of 
the well-known Jacobian conjecture (see \cite{K, BCW, E}). 
See also \cite{DEZ, EKC}. However, a more interesting aspect 
of the notion is that it provides a natural but 
highly non-trivial generalization of the notion 
of ideals. Currently, this new notion has not 
been studied (nor understood) for the most of rings including 
the most of finite rings and finite dimensional algebras 
over a field. 

\begin{defi} \cite[p.\,247]{MS} \label{Def-Rad}
Let $V$ be an $R$-subspace (or a subset) of an $R$-algebra $\cA$. 
We define the {\it radical} $\rad(V)$ of $V$ to be 
\begin{align}
\rad(V)\!:=\{ a\in \cA\,|\, a^m\in V \text{ for all } m\gg0\}.
\end{align}
\end{defi}

When $\cA$ is commutative and $V$ is an ideal of $\cA$, $\rad(V)$ 
coincides with the radical of the ideal $V$, which is defined as 
$\rad(V)=\{ a\in \cA\,|\, a^m\in V \text{ for some } m\ge 1\}$. 
So this new notion generalizes the radical of ideals and is 
interesting on its own right. It is also crucial for the study of MSs. 
For example, it is easy to see that {\it every $R$-subspace 
$V$ of an $R$-algebra $\cA$ with $\rad(V)\subseteq \nil(\cA)$
(equivalently, $\rad(V)=\nil(\cA)$, since   
$\nil(\cA)\subseteq \rad(V)$ by definition) is a MS of $\cA$, 
where $\nil(\cA)$ denotes the set of  all nilpotent elements of $\cA$.}
We will frequently use this fact (implicitly) throughout this paper.

%
 
Recent studies show that many MSs arise from 
the images of differential operators, especially, 
from the images of locally finite or locally 
nilpotent derivations, of certain associative 
algebras (e.g., see \cite{IC, GIC, EWZ, EZ1, EZ2} 
and \cite{Open-LFNED}--\cite{OneVariableCase}, etc.). 
Then one natural question is the following:
\begin{oprob} 
For which differential operator $L$ of $\cA$,  
the kernel $\Ker L$ of $L$ is a MS of $\cA$?
\end{oprob}
Note that differential operators are among the most classical and fundamental  subjects in mathematics. They have been extensively studied not only in theories of ODE and PDE, $D$-modules, differential or complex manifolds, etc., but also in many other different areas such as general theories of rings and algebras (e.g., see \cite{Kharchenko1}, \cite{Kharchenko2} and the references therein). Nevertheless, it seems that the question above and the radical of the kernel of differential operators have not been studied before! It is presumably because MSs and the radical of (arbitrary) subspaces are still relatively very new notions. After all, they were introduced in \cite{GIC, MS} only about a decade ago.    


In this paper we study the open problem above. More precisely, 
we study the radicals of the kernels of some (ordinary or partial) differential operators of some $R$-algebras $\cA$, and show that for 
certain differential operators $\Phi$ of $\cA$,  
the kernel $\Ker \Phi$ is indeed a MS of $\cA$. 
We also apply some results 
proved in this paper to study $R$-derivations of $\cA$ 
that are locally algebraic or locally integral  over $R$  
(see Definition \ref{Def-LA}). In particular, 
we show that if $R$ is an integral domain of characteristic zero and 
$\cA$ is reduced and torsion-free as an $R$-module, then 
$\cA$ has no nonzero $R$-derivation that is locally algebraic over $R$ 
(see Theorem \ref{NoncommDomainCase-1}).  
Finally, we also show a formula for the determinant of a  
differential vandemonde matrix over commutative algebras 
(see Proposition \ref{DiffVandeMonde}). This formula not only 
provides some information for the radicals of the kernels 
of ordinary differential operators of commutative algebras, 
but also is interesting on its own right. \\


%
%
%
%
%
%

{\bf Arrangement and Content}: In Section \ref{S2}, we assume that $\cA$ is commutative 
and derive some necessarily conditions 
for the elements in the radical of the kernel of an arbitrary 
differential operator of $\cA$ (see Theorem \ref{CommCase} and Corollary \ref{Corol-2.4}). 
Consequently, for every differential operator $\Phi\in\cW(\cA)$ such that $\Phi 1_\cA$ is not zero nor a zero-divisor of $\cA$, the kernel $\Ker \Phi$ is indeed a MS of $\cA$.  
  
In Section \ref{S3}, we drop the commutativity assumption on $\cA$ but assume that $(R, +)$ is torsion-free and $\cA$ is reduced and torsion-free as an $R$-module. 
We first show in Theorem \ref{Thm3.1} some necessary conditions satisfied by the elements in the radical of the kernel of a differential operator $P(D)$ of $\cA$, where $P(\cdot)$ is a polynomial over $R$ in $n$ commutative free variables and $D=(D_1, D_2, \dots, D_n)$ are $n$ commuting $R$-derivations of $\cA$ such that for each $1\le i\le n$, $\cA$ can be decomposed as a direct sum of the generalized 
eigen-subspaces of $D_i$.   

We then show in Proposition \ref{NonCommLF-Case} that if $R$ is an integral domain of characteristic zero, then the conclusions in Theorem \ref{Thm3.1} hold also for the differential operators of $\cA$ which are multivariate polynomials over $R$ 
in commuting locally finite $R$-derivations of $\cA$. Finally, we show in 
Proposition \ref{NonCommOneVariaCase} that the similar conclusions as those in Proposition \ref{NonCommLF-Case} (with the same assumptions on $R$ and $\cA$) hold also for all ordinary differential operators of $\cA$.
Consequently, for all the differential operators $\Phi$ in Theorem \ref{Thm3.1} and Propositions \ref{NonCommLF-Case}, \ref{NonCommOneVariaCase} with $\Phi 1_\cA\ne 0$,  $\Ker \Phi$ is indeed a MS of 
$\cA$.  

In Section \ref{S4}, we apply some results proved in Sections \ref{S2} 
and \ref{S3} to study some properties of $R$-derivation of $\cA$, which  are locally algebraic or locally integral over $R$ (see Definition \ref{Def-LA}). We first show in Theorem \ref{CommLocIntCase} that if $\cA$ is commutative and $(\cA, +)$ is torsion-free, then every locally integral $D$ of $\cA$ has its image in the nil-radical $\nil(\cA)$ of $\cA$. We then show in Theorem \ref{NoncommDomainCase-1} that if $R$ is an integral domain of characteristic zero and $\cA$ is reduced and torsion-free 
as an $R$-module (but not necessarily commutative), then $\cA$ has no nonzero $R$-derivation that is locally algebraic over $R$. 

In Section \ref{S5}, we assume that $\cA$ is commutative and first 
show in Proposition \ref{DiffVandeMonde} a formula for the determinant of a differential vandemonde matrix over $\cA$. We then apply this formula 
in Proposition \ref{Propo-5.4} to derive more necessary conditions satisfied by the elements in the radicals of the kernels of 
an ordinary differential operator of $\cA$.   
we point out in Remark \ref{OtherApp} that the formula derived in Proposition \ref{DiffVandeMonde} can also be used to derive formulas for the determinants of several other families of matrices. Therefore the formula is also interesting on its own right.

\renewcommand{\theequation}{\thesection.\arabic{equation}}
\renewcommand{\therema}{\thesection.\arabic{rema}}
\setcounter{equation}{0}
\setcounter{rema}{0}

\section{\bf The Commutative Algebra Case}\label{S2}

In this section, unless stated otherwise, {\it $R$ denotes a unital commutative ring, $\cA$ a commutative unital $R$-algebra and $\xi=(\xi_1, \xi_2, \dots, \xi_n)$ $n$ noncommutative free variables}. 
We denote by $\cA\langle\xi\rangle$ the (noncommutative) polynomial algebra in $\xi$ over $\cA$, and by $\p_i$ $(1\le i\le n)$ the $\cA$-derivation $\p/\p \xi_i$ of $\cA\langle\xi\rangle$.  

Once and for all, we fix in this section a nonzero $P(\xi)\in\cA\langle\xi\rangle$
and $n$ $R$-derivations $D_i$ $(1\le i\le n)$ of $\cA$. 
Write $D=(D_1, D_2, \dots, D_n)$ and 
$P(\xi)=a_0+\sum_{k=1}^d P_k(\xi)$ for some $a_0\in\cA$, $d\ge 1$ 
and homogeneous polynomials $P_k(\xi)$ ($1\le k\le d$) 
of degree $k$ in $\xi$. 

For each $u\in \cA$, we set 
$\nabla_D u\!:=(D_1u, D_2u, \dots, D_nu)$, 
and call it the {\it gradient} of $u$ 
with respect to $D$. When $D$ is clear in the context, 
we will simply write $\nabla_D u$ as $\nabla u$.  

We define $P(D)$ and $P(\nabla u)$ 
by first {\it writing $P(\xi)$ with all the coefficients of $P(\xi)$ 
on the most left of the monomials in $\xi$}, and then replacing $\xi_i$ 
by $D_i$ and $D_iu$, respectively, for each $1\le i\le n$.   

Note that every differential operator $\Psi$ 
in the Weyl algebra $\cW_R(\cA)$ of $\cA$ can be written as 
$\Psi=P(D)$ for some $R$-derivations $D_i$ $(1\le i\le n)$ 
of $\cA$ and $P(\xi)\in \cA\langle\xi\rangle$.

The main result of this section is the following:  

\begin{theo}\label{CommCase}
With the setting as above, let $u\in \cA$ be such that $u^m\in \Ker P(D)$ 
for all $1\le m\le d$. Then  
\begin{align} 
a_0 u^d=(-1)^d d!\,P_d(\nabla u).\label{CommCase-eq1} 
\end{align}  

Furthermore, if $u^{d+1}$ also lies in $\Ker P(D)$, 
then  
\begin{align}\label{CommCase-eq2}
a_0 u^{d+1}=0.  
\end{align}  
\end{theo}

To show the theorem above, we need first the following two lemmas. 
The first lemma  can be easily verified by using 
the mathematical induction, 
which is similar as the proof for the usual 
binomial formula. So we here skip its proof. 

\begin{lemma}\label{CommCase-Lma1}
Given $u\in \cA$, let $\ell_u:\cA \to \cA$ be the map such that 
$\ell_u(a)=au$ for all $a\in \cA$. 
Define $\ad_u:\cW (\cA)\to \cW (\cA)$ by setting 
$\ad_u(\Psi)=[\ell_u, \Psi]\!:=\ell_u\Psi-\Psi \ell_u$ for all 
$\Psi \in \cW (\cA)$. 
Then for all $\Phi\in\cW (\cA)$ and $k\ge 1$, we have 
\begin{align}\label{CommCase-Lma1-eq1}
(\ad_{u})^k (\Phi)=\sum_{i=0}^k (-1)^i \binom{k}i u^{k-i} \Phi \circ \ell_u^i. 
\end{align}
\end{lemma}

\begin{lemma}\label{CommCase-Lma2}
Let $u\in \cA$ and $d=\deg P(\xi)$.   
Then the following statements hold:
\begin{enumerate}
\item[$1)$] there exists $Q(\xi)\in \cA\langle\xi\rangle$ with either $Q(\xi)=0$ or 
$\deg Q(\xi)\le d-2$ such that  
\begin{align}
\ad_{-u} P(D)=\sum_{i=1}^n (D_i u)(\p_i P)(D)+ Q(D). 
\end{align}

\item[$2)$]\hspace{25mm} 
$
(\ad_{-u})^d P(D)=d!P(\nabla u).  
$
\end{enumerate} 
\end{lemma}

\pf $1)$ First, if $d=\deg P(\xi)=0$, then the statement holds trivially, 
since $\cA$ is commutative and hence $\ad_{-u} P(\xi)=0$. So 
we assume $d=\deg P(\xi)\ge 1$. By the linearity of $\ad_u$ and $D_i$'s 
and also by the commutativity of $\cA$ we may assume 
$P(\xi)=\xi_{i_1}\xi_{i_2}\cdots\xi_{i_d}$ with $1\le i_j\le n$  
for all $1\le j\le d$. 

We use the induction on $d\ge 1$. 
If $d=1$, then $\ad_{-u}D_{i_1}=\ell_{D_{i_1}u}$. Hence the statement holds by 
choosing $Q(\xi)=0$. Assume that the statement holds for all $1\le d \le m-1$ and 
consider the case $d=m$. 

Since $\ad_{-u}$ is a derivation of $\cW(\cA)$, we have    
\begin{align*}
&\ad_{-u}P(D)=\sum_{j=1}^m  D_{i_1} \cdots (\ad_{-u}{D_{i_j}}) \cdots D_{i_m}
=  \sum_{j=1}^m  D_{i_1} \cdots (\ell_{D_{i_j}u}) \cdots D_{i_m}\\
&=\sum_{j=1}^m  (\ell_{D_{i_j}u}) D_{i_1} \cdots \widehat{D_{i_j}} \cdots D_{i_m}+
\sum_{j=2}^m  [D_{i_1} \cdots D_{i_{j-1}}, \, \ell_{D_{i_j}u}] D_{i_{j+1}} \cdots D_{i_m},
\end{align*}
where $\widehat{D_{i_j}}$ means that the term $D_{i_j}$ is omitted. Thus 
\begin{align*}
\ad_{-u}P(D)=\sum_{j=1}^m  (D_{i_j}u) & D_{i_1} \cdots \widehat{D_{i_j}} \cdots D_{i_m}\\  & +
\sum_{j=2}^m  (\ad_{-D_{i_j}u}(D_{i_1} \cdots D_{i_{j-1}})) D_{i_{j+1}} \cdots D_{i_m}.
\end{align*}
Applying the induction assumption to the terms 
$\ad_{-D_{i_j}u}(D_{i_1} \cdots D_{i_{j-1}})$ 
$(2\le j\le m)$ in the sum above we see that there exists 
$Q(\xi)\in \cA\langle\xi\rangle$ with $Q(\xi)=0$ or $\deg Q(\xi)\le m-2$ 
such that 
\begin{align*}
\ad_{-u}P(D)&=\sum_{j=1}^m  (D_{i_j}u) D_{i_1} \cdots \widehat{D_{i_j}} \cdots D_{i_m}+Q(D)\\
&=\sum_{i=1}^n  (D_i u) (\p_i P)(D)+Q(D).
\end{align*}
Hence by the induction statement $1)$ follows.

$2)$ First, by statement $1)$ it is easy to see that 
$(\ad_{-u})^d P(D)=(\ad_{-u})^d P_d(D)$. Then  
by the linearity of $\ad_u$ and $D_i$'s 
and also by the commutativity of $\cA$ 
we may assume $P_d(\xi)=\xi_{i_1}\xi_{i_2}\cdots\xi_{i_d}$ 
with $1\le i_j\le n$  
for all $1\le j\le d$. Applying statement $1)$ ($d$ times) 
we have 
\begin{align*}
(\ad_{-u})^d P(D)&=
\sum_{1\le k_1, k_2, \dots, k_d\le n}  
(D_{k_1}u)(D_{k_2}u)\cdots (D_{k_d}u)
(\p_{k_1}\p_{k_2}\cdots \p_{k_d} P).
\end{align*}
Then by the equation above and the 
commutativity of $\cA$, it is easy to see that the equation in statement 
$2)$ follows. 
\epfv

Now we can prove the main result of this section. \\

\underline{\bf Proof of Theorem \ref{CommCase}:}\, 
By Eq.\,(\ref{CommCase-Lma1-eq1}) and Lemma \ref{CommCase-Lma2}, $2)$ we have 
\begin{align}\label{CommCase-peq1}
d!P_d(\nabla u)=(-1)^d\sum_{i=0}^d (-1)^i \binom{d}i u^{d-i} P(D) \circ \ell_u^i.
\end{align}
By applying both sides of the equation above to $1_\cA\in \cA$ and then 
using the condition $u^i\in \Ker P(D)$ $(1\le i\le d)$, we get $d!P_d(\nabla u)=(-1)^d  u^d P(D)\cdot 1_\cA$. It is well-known and also easy to check that every derivation of a commutative ring  annihilates the identity element of the ring.  Hence $d!P_d(\nabla u)= (-1)^d  u^da_0$, 
i.e., Eq.\,(\ref{CommCase-eq1}) follows. Similarly, by applying Eq.\,(\ref{CommCase-peq1}) above 
to $u\in \cA$ and using the condition $u^i\in \Ker P(D)$ $(1\le i\le d+1)$, we get 
$d!P_d(\nabla u)u =0$. Then by  
Eq.\,(\ref{CommCase-eq1}) we get Eq.\,(\ref{CommCase-eq2}).
\epfv

Two consequences of Theorem \ref{CommCase} are as follows.

\begin{corol}\label{New-Corol-2.4}
Let $D$, $P(\xi)$ and $a_0$ be as in Theorem \ref{CommCase}. If $P(D)=0$, then $a_0=0$ 
and $d!\,P_d(\nabla u)=0$ for all $u\in\cA$.
\end{corol}
\pf Since $P(D)=0$, we have $\cA=\Ker P(D)=\rad(\Ker P(D))$. Applying Eq.\,(\ref{CommCase-eq2}) to 
$u=1$ we get $a_0=0$. Then the corollary follows immediately from Eq.\,(\ref{CommCase-eq1}).
\epfv

\begin{corol}\label{Corol-2.4}
Let $D$, $P(\xi)$, $a_0$ be as in Theorem \ref{CommCase}, and $\nil(\cA)$ 
the nil-radical of $\cA$, i.e., the set of all nilpotent elements of $\cA$.     
Then the following statements hold:
\begin{enumerate}
  \item[$1)$] $\rad(\Ker P(D))\subseteq \rad(\ann (a_0))$, where $\ann(a_0)$ is the set of the elements $b\in \cA$ such that $a_0b=0$;
\item[$2)$] if $a_0$ is not zero nor a zero-divisor of $\cA$, then 
$\rad(\Ker P(D))=\nil(\cA)$ and $\Ker P(D)$ is a MS of $\cA$;   
\item[$3)$] if $a_0=0$, then we have
\begin{align*}
\rad(\Ker P(D))\subseteq \rad\big(\{ u\in\cA\,|\, d!P_d(\nabla u)=0\}\big). 
\end{align*}
In particular, 
if $n=1$, i.e., $D$ is a single derivation of $\cA$, and 
the leading coefficient of $P(\xi)$ is not a zero-divisor of $\cA$, 
then 
\begin{align*}
\rad(\Ker P(D))\subseteq \rad\big(\{ u\in\cA\,|\, d! Du\in \nil (A)\}\big). 
\end{align*}
\end{enumerate}
\end{corol}

\pf $1)$ Let $u\in \rad(\Ker P(D))$. Then there exists $N\in\bN$ 
such that $u^m \in \Ker P(D)$ for all $m\ge N$. 
In particular, $u^{Nk}=(u^N)^k\in \Ker P(D)$ for all $k\ge 1$. 
Then by Theorem \ref{CommCase} we have $a_0 u^N=0$, and hence 
$a_0u^\ell=0$ for all $\ell\ge N$. Thus $u\in \rad(\ann (a_0))$, 
and statement $1)$ follows.

$2)$ Since $a_0$ is not zero nor a zero-divisor of $\cA$, 
we have $\ann (a_0)=\nil(\cA)$. 
By Definition \ref{Def-Rad} it is easy to see that 
$\rad( \nil(\cA) )\subseteq \nil(\cA)$, and $\nil(\cA)\subseteq \rad(V)$ 
for all $R$-subspace $V\subseteq \cA$. Then by statement $1)$ 
we have $\rad(\Ker P(D))=\rad(\nil(\cA) )=\nil(\cA)$ and 
$\Ker P(D)$ is a MS of $\cA$.

$3)$ Let $v\in \rad(\Ker P(D))$. Then there exists $N\in\bN$ 
such that $v^m \in \Ker P(D)$ for all $m\ge N$. 
In particular, for all $m\ge N$ we have 
$v^{mk}=(v^m)^k\in \Ker P(D)$ for all $k\ge 1$.
Then by Eq.\,(\ref{CommCase-eq1}) and the condition $a_0=0$ we have 
$v^m\in \{ u\in\cA\,|\, d!P_d(\nabla u)=0\}$ (for all $m\ge N$). Thus 
$v\in \rad \big(\{ u\in\cA\,|\, d!P_d(\nabla u)=0\}\big)$. 

Furthermore, assume that $n=1$ and the leading coefficient of $P(\xi)$ 
is not a zero-divisor of $\cA$. Then for all $u\in \cA$ with $d!P_d(\nabla u)=0$ 
we have $d!(Du)^d=0$. Hence $(d!Du)^d=0$ and $d!Du \in \nil(\cA)$. 
Then in this case we have 
$v\in \rad \big(\{ u\in\cA\,|\, d! Du\in \nil (A) \}\big)$. 
\epfv

\begin{exam} 
Let $R=\bC$ and $\cA$ the $\bC$-algebra of all smooth complex valued functions $f(x)$ 
over $\bR$. Let $D=\frac{d}{dx}$. Then for each nonzero univariate 
polynomial $P(\xi)\in \bC[\xi]$,  $\Ker P(D)$ is the set 
of solutions $f(x)\in \cA$ of the ordinary differential 
equation $P(D)f=0$. 

Let $\lambda_i$ $(1 \le i \le k)$ be the set of all distinct roots 
of $P(\xi)$ in $\bC$ with multiplicity $m_i$. Then it is well-known 
in the theory of ODE (e.g., see \cite{L} or any other standard 
textbook on ODE) that $\Ker P(D)$ is the $\bC$-subspace of $\cA$ 
spanned by $x^je^{\lambda_ix}$ for all $1\le i\le k$ and 
$1\le j\le m_i$.   

Then it is readily verified directly from the fact above that 
$\rad(\Ker P(D))=\{0\}$ if $P(0)\ne 0$; and 
$\rad(\Ker P(D))=\bC$ if $P(0)=0$. Consequently, 
Theorem \ref{CommCase}, Corollary \ref{Corol-2.4} and also  
Proposition \ref{Propo-5.4} in Section \ref{S5} all 
hold in this case.
\end{exam}

Furthermore, from the example below we get a family of examples of MSs  
from the solution spaces of linear Partial Differential Equations.

\begin{exam} \label{PDE-MS}
Let $R=\bC$ or $\bR$, and $\cA$ be the $R$-algebra of all smooth $R$-valued functions 
$f(x)$ over an open subset of $\bR^n$ (or let $\cA$ be the polynomial algebra 
in $n$ commutative free variables over $R$). Let $D_i=\p/{\p x_i}$ $(1\le i\le n)$ 
and $D=(D_1, D_2, ..., D_n)$. Then for every partial differential operator $\Phi$ 
of $\cA$, there exists a polynomial $P(\xi)\in \cA[\xi]$ 
such that $\Phi=P(D)$. Then by Corollary \ref{Corol-2.4}, 
$2)$ we see that $\Ker \Phi$, or equivalently, the solution space in $\cA$ of the 
PDE: $\Phi f=0$, is a MS of $\cA$ as long as $\Phi 1_{\cA}\ne 0$.
\end{exam}

We end this section with the following two remarks.

First, we will show in Propositions \ref{NonCommOneVariaCase} and \ref{Propo-5.4} that 
for the ordinary differential operators $\Psi$ of certain $R$-algebras $\cA$ 
(not necessarily commutative),  
the radical $\rad(\Ker \Psi)$ also satisfies some other necessarily conditions 
(other than those in Theorem \ref{CommCase} and Corollary \ref{Corol-2.4}).   

Second, Theorem \ref{CommCase} and Corollary \ref{Corol-2.4} cannot be generalized to  differential operators of noncommutative algebras, which can be seen 
from the following:  

\begin{exam} 
Let $X$, $Y$ be two noncommutative free variables and $R\langle X, Y\rangle$ the polynomial algebra in $X$ and $Y$ over $R$. Let $J$ be the two-sided ideal of $R\langle X, Y\rangle$ generated by $Y^2$ and 
$\cA\!:=R\langle X, Y\rangle/J$. 
Let $D=\p/\p X$ and $P(\xi)=1-X \xi\in\cA\langle\xi\rangle$. Then $P(D)=\I-\ell_X D$, 
where $\I$ denotes the identity map of $\cA$, and $\ell_X$ the multiplication map by $X$ from the left. Let $f=XY\in \cA$. Then it is readily checked that for all $m\ge 1$, we have  
$P(D)(f^m)=0$ but $P(D)(Xf^m)=-Xf^m\neq 0$. 
Therefore, $0\ne f\in \rad(\Ker P(D))$ and 
$\Ker P(D)$ is not a MS of $\cA$. 
\end{exam}

%
%
%
%
%

\section{\bf Some Cases for Non-Commutative Algebras}\label{S3}

In this section, unless stated otherwise, {\it $R$ denotes a commutative ring such that the abelian group $(R, +)$ is torsion-free, and $\cA$ an $R$-algebra (not necessarily commutative) 
that is torsion-free as an $R$-module.}  

We denote by $\I_\cA$, or simply $\I$,  
the identity map of $\cA$, and by $\nil(\cA)$ the set of all 
nilpotent elements of $\cA$. We say $\cA$ is {\it reduced} 
if $\nil(\cA)=\{0\}$.
Furthermore, for each $a\in \cA$, we denote by $\ann_\ell(a)$  
 the set of elements $b\in \cA$ such that $ab=0$. 

Let $D$ be an $R$-derivation of $\cA$. We say that 
$\cA$ is {\it decomposable w.r.t. (with respect to) $D$} if $\cA$ can be written as a 
direct sum of the generalized eigen-subspaces of $D$. 
More precisely, $\cA=\oplus_{\lambda\in H}\cA_{\lambda}$, where 
$H$ is the set of all generalized eigenvalues of $D$ in $R$ 
and $\cA_\lambda=\sum_{i=1}^\infty \Ker (D-\lambda\I)^i$ 
for each $\lambda\in H$. 
It is easy to verify inductively that for all $m\ge 1$, 
$a, b\in \cA$ and $\lambda, \mu\in R$, we have 
\begin{align}\label{LamMu-Rule}
\big(D-(\lambda+\mu)\I\big)^m(ab)=\sum_{i=0}^m\binom{m}{i} 
\big((D-\lambda\I)^ia\big)\big((D-\mu\I)^{m-i}b\big).
\end{align}
Then by the identity above we have that 
$\cA_{\lambda}\cA_\mu\subseteq \cA_{\lambda+\mu}$  for all 
$\lambda, \mu\in H$. 
In other words, the decomposition $\cA=\oplus_{\lambda\in H}\cA_{\lambda}$  
is actually an additive $R$-algebra grading of $\cA$. 

Some examples of $R$-derivations with respect to which $\cA$ is decomposable are 
{\it semi-simple} $R$-derivations, for which $\cA_\lambda$ ($\lambda\in H$) coincides 
with the eigenspace of $D$ corresponding to the eigenvalue $\lambda$ of $D$, 
and also {\it locally finite} derivations when the base ring $R$ is 
an algebraically closed field (e.g., see \cite[Proposition $1.3.8$]{E}) .   

Throughout this section $D_i$ $(1\le i\le n)$ stand for $n$ commuting $R$-derivations of $\cA$, 
i.e.,  $D_iD_j=D_jD_i$ for all $1\le i, j\le n$, 
such that $\cA$ is decomposable w.r.t. each $D_i$. 
Then there exists a semi-subgroup $\Lambda$ of the abelian group $(R^n, +)$ 
such that 
\begin{align}\label{Decom-A}
\cA=\oplus_{\lambda\in \Lambda}\cA_{\lambda},
\end{align}
where for each $\lambda=(k_1, k_2, \dots, k_n)\in\Lambda$,  
\begin{align}
\cA_{\lambda}=\bigcap_{i=1}^n \big(\sum_{j=1}^{\infty} 
\Ker (D_i-k_i\I)^j\big).
\end{align}
In particular,  
\begin{align}
\cA_0=\bigcap_{i=1}^n \big(\sum_{j=1}^{\infty} 
\Ker D_i^j\big).
\end{align}

Note also that each $\cA_\lambda$ $(\lambda\in \Lambda)$ is invariant under $D_i$ 
$(1\le i\le n)$, and $\cA_\lambda\cA_\mu \subseteq\cA_{\lambda+\mu}$ for all 
$\lambda, \mu\in \Lambda$.

Now, let $\xi=(\xi_1, \xi_2, \dots, \xi_n)$ be $n$ 
commutative free variables and $R[\xi]$ the polynomial algebra in $\xi$ 
over $R$. We set $D\!:=(D_1, D_2, \dots, D_n)$ and fix a polynomial 
$0\ne P(\xi)\in R[\xi]$. 

Write $P(\xi)=\sum_{k=0}^d P_k(\xi)$ for some $d\ge 0$ 
and homogeneous polynomials $P_k(\xi)$ ($1\le k\le d$) 
of degree $k$ in $\xi$. Let $P(D)$ be the differential operator of $\cA$ 
obtained by replacing $\xi_i$ by $D_i$ $(1\le i\le n)$. 
Since $D_i$'s are $R$-derivations and commute with one anther, 
$P(D)$ is well-defined. 

The first main result of this section is the following theorem which in some sense 
extends Theorem \ref{CommCase} to the differential operator $P(D)$ of 
the $R$-algebra $\cA$ (that is not necessarily commutative). 

\begin{theo}\label{Thm3.1} 
With the setting as above, assume further that $\cA$ is reduced. 
Then the following statements hold:
\begin{enumerate}
  \item[$1)$] if $P(0)= 0$ and  $P_k(\xi)$ $(1\le k\le d)$ have no nonzero common zeros in $R^n$, then  
$ 
\rad(\Ker P(D))\subseteq \cA_0;
$ 
\item[$2)$] if $P(0)\ne 0$, then $\rad(\Ker P(D))=\{0\}$, 
and $\Ker P(D)$ is a MS of $\cA$.
\end{enumerate}
\end{theo}

In order to prove the theorem above, we need first to show some lemmas.  

\begin{lemma}\label{Lma3.1}
Let $R$ be an arbitrary commutative ring and 
$\cA$ an $R$-algebra that is torsion-free as an $R$-module. 
Let $D$ and $P(\xi)$ be fixed as above. 
Then the following statements hold: 
\begin{enumerate}
  \item[$1)$] $\Ker P(D)$ is homogeneous w.r.t. the grading of $\cA$ in Eq.\,(\ref{Decom-A}), i.e., 
\begin{align}\label{Lma3.1-eq1}
\Ker P(D)=\bigoplus_{\lambda\in \Lambda} (\cA_{\lambda}\cap \Ker P(D)).
\end{align}
\item[$2)$] Let $\cZ_\Lambda(P)$ be the set of $\lambda\in \Lambda$ such that 
$P(\lambda)=0$. Then 
\begin{align}\label{Lma3.1-eq2}
\Ker P(D)\subseteq \bigoplus_{\lambda\in \cZ_\Lambda(P)}\cA_\lambda.
\end{align} 
\end{enumerate}
\end{lemma}  

\pf $1)$ Since for each $\lambda\in \Lambda$, 
$\cA_\lambda$ is preserved by $D_i$ $(1\le i\le n)$, and hence 
is also preserved by $P(D)$, from which Eq.\,(\ref{Lma3.1-eq1}) 
follows. 
 
$2)$ Let $0\ne u\in\Ker P(D)$ and write 
$u=\sum_{i=1}^\ell u_{\lambda_i}$ for some distinct 
$\lambda_i\in \Lambda$   
and  $u_{\lambda_i}\in\cA_{\lambda_i}$ ($1\le i\le \ell$). 
Then by Eq.\,(\ref{Lma3.1-eq1})   
we have that $u_{\lambda_i}\in \Ker P(D)$ for all $1\le i\le \ell$.  
So we may assume $\ell=1$ and $u\in \cA_{\lambda}$ 
for some $\lambda\in \Lambda$. 

Write $\lambda=(k_1, k_2, \dots, k_n)$. 
For each $1\le j\le n$, we define a non-negative integer 
$r_j$ as follows. 
First, let $r_n$ be the greatest  
non-negative integer such that $(D_n-k_n\I)^{r_n}u \ne 0$, and inductively, for each $1\le j\le n-1$, let $r_j$ be the greatest non-negative  
integer such that 
$$(D_j-k_j\I)^{r_j}\left(\prod_{s=j+1}^n (D_s-k_s\I)^{r_s}\right) u\ne 0.$$
Set $\tilde u\!:=\left(\prod_{j=1}^n(D_j-k_j\I)^{r_j}\right)u$. 
Then $0\ne \tilde u\in \cA_{\lambda}$, $\tilde u\in \Ker P(D)$, and $D_j\tilde u=k_j\tilde u$ for all $1\le j\le n$. Hence $0=P(D)\tilde u=P(\lambda)\tilde u$. Since $\cA$ is torsion-free as an $R$-module, 
we have  $P(\lambda)=0$, as desired. 
\epfv

\begin{defi}\label{Def-ExtremeElt}
Let $B$ be a subset of $R^n$ and $\lambda\in A$. We say $\lambda$  
is {\it an extremal element} of $B$ if for all $m\ge 1$, $m\lambda$ can not be written as a linear combination of other elements of $B$  
with positive integer coefficients 
whose sum is less than or equal 
to $m$. 
\end{defi}

The following lemma should be known. But for the sake of completeness, 
we include here a direct proof.

\begin{lemma}\label{Lma3.2}
Let $R$ be a commutative ring such that the abelian group $(R, +)$ 
is torsion free. Then 
every nonempty finite subset $B$ of $R^n$ has at least one 
extremal element.
\end{lemma}

\pf Write $B=\{\lambda_1, \lambda_2,\dots, \lambda_n\}$ with $\lambda_i\ne\lambda_j$ for all 
$1\le i\ne j\le n$. We use induction on $n$. 
If $n=1$, there is nothing to show. 
So we assume $n\ge 2$. 

Consider first the case $n=2$ with $\lambda_2\ne 0$.
If the lemma fails, then $m_1\lambda_1=k_1\lambda_2$ and 
$m_2\lambda_2=k_2\lambda_1$ 
for some $m_i, k_i\ge 1$ with $k_i\le m_i$. Then 
$m_1 m_2\lambda_2=m_1(k_2\lambda_1)=
k_2(m_1\lambda_1)=k_1k_2\lambda_2$. Hence $m_1m_2=k_1k_2$, 
for $\lambda_2\ne 0$ and $(R, +)$ is torsion-free, from which we have $m_1=k_1$ 
(and $m_2=k_2$). By the assumption that $(R, +)$ is torsion-free again, 
we have $\lambda_1=\lambda_2$. Contradiction.

Now assume the lemma holds for all $2\le n\le k$ and consider 
the case $n=k+1$. If $\lambda_{k+1}$ is an extremal point of $A$, 
then there is nothing to show. Assume otherwise. Then there exist 
$m \ge 1$ and  $c_i\in \bN$ ($1\le i\le k$) such that 
\begin{align}
m \lambda_{k+1}=\sum_{i=1}^k c_i\lambda_i, \label{Lma3.2-peq1}  \\ 
1\le \sum_{i=1}^k c_i\le m. \label{Lma3.2-peq2}
\end{align}

By the induction assumption the set 
$B'\!:=\{\lambda_1, \lambda_2,\dots, \lambda_k\}$ has an extremal element, 
say, $\lambda_1$. 
We claim that $\lambda_1$ is also an extremal point of the set $B$. Otherwise, 
there exist $q\ge1$ and $c_j'\in \bN$ ($2\le j\le k+1$) such that 
\begin{align}
q&\lambda_1=c_{k+1}'\lambda_{k+1}+\sum_{j=2}^{k}c_j' \lambda_j, \label{Lma3.2-peq3}\\
&\,\,1\le c_{k+1}'+\sum_{j=2}^{k} c_j'\le q. \label{Lma3.2-peq4}
\end{align}
Then by Eqs.\,(\ref{Lma3.2-peq3}) and (\ref{Lma3.2-peq1}) we have 
\begin{align}
m q \lambda_1&=m c'_{k+1}\lambda_{k+1}+
m \sum_{j=2}^{k} c_j'\lambda_j \label{Lma3.2-peq5} \\
 &=c'_{k+1} \sum_{i=1}^k c_i\lambda_i +
m \sum_{j=2}^{k} c_j'\lambda_j.  \nno 
\end{align}
For the sum of all the coefficients of the linear combination on the right hand side of the equation above, by Eqs.\,(\ref{Lma3.2-peq2}) and (\ref{Lma3.2-peq4}) we have  
\begin{align}
1&\le c'_{k+1} \sum_{i=1}^k c_i +
m \sum_{j=2}^k c_j' \le c'_{k+1}m +
m \sum_{j=1}^k c_j'   \nno  \\&
= m (c'_{k+1}+\sum_{j=2}^k c_j')\le m q.  \label{Lma3.2-peq6}
\end{align}
Then by Eqs.\,(\ref{Lma3.2-peq5}) and (\ref{Lma3.2-peq6}), 
$\lambda_1$ is not an extremal element of $B'$, which contradicts the choice 
of $\lambda_1$. Therefore $\lambda_1$ is an extremal point of $B$, 
and the lemma follows. 
\epfv

\begin{lemma}\label{Lma3.3}
Let $0\ne u\in \rad(\Ker P(D))$ and write $u=\sum_{i=1}^\ell 
u_{\lambda_i}$ for some distinct $\lambda_i\in \Lambda$ $(1\le i\le \ell)$ and 
$0\ne u_{\lambda_i}\in \cA_{\lambda_i}$. 
 Then for each  
 extremal element $\lambda_j$ of the set 
$\{\lambda_i\, |\, 1\le i\le \ell\}$, 
either $u_{\lambda_j}$ is nilpotent, or 
$P_k(\lambda_j)=0$ for all $0\le k\le d$. 
\end{lemma}

\pf Assume that $u_{\lambda_j}$ is not nilpotent. Since 
$\lambda_j$ is an extremal element of the set $\{\lambda_i\, |\, 1\le i\le \ell\}$, 
it is easy to see that for each $m\ge 1$, the homogeneous component of 
$u^m$ in $\cA_{m\lambda_j}$ is equal to $u_{\lambda_j}^m$. 
Since $u^m\in \Ker P(D)$ when $m\gg 0$, 
by Lemma \ref{Lma3.1}, $1)$ and $2)$ 
we have $u_{\lambda_j}^m\in \Ker P(D)$ and 
$P(m\lambda_j)=0$ for all $m\gg 0$. 
More explicitly, for all $m\gg 0$, we have 
$$
0=P(m\lambda_j)=\sum_{k=0}^d m^k P_k(\lambda_j).
$$
Since $(R, +)$ is torsion-free, by the vandemonde determinant 
we have $P_k(\lambda_j)=0$ for all $0\le k\le d$. 
\epfv
 
\underline{\bf Proof of Theorem \ref{Thm3.1}}:  
Let $0\ne u\in \rad(\Ker P(D))$ and write $u=\sum_{i=1}^\ell 
u_{\lambda_i}$ for some distinct $\lambda_i\in \Lambda$ $(1\le i\le \ell)$ and 
$0\ne u_{\lambda_i}\in \cA_{\lambda_i}$. 
Let $B$ be the set of all nonzero $\lambda_i$ 
($1\le i\le \ell$). 

If $B\ne \emptyset$, then by Lemma \ref{Lma3.2}, $B$ has at least one extremal element, say $\lambda_j$. Then by Definition \ref{Def-ExtremeElt}, 
$\lambda_j$ is also an extremal element of the set $B\cup\{0\}$. 
Since $\cA$ is reduced, $u_{\lambda_j}$ is not nilpotent. 
Then by Lemma \ref{Lma3.3}, $P_k(\lambda_j)=0$ 
for all $0\le k\le d$.

If $P(0)=0$, and $P_k(\xi)$ $(1\le k\le d)$ have no nonzero 
common zero in $R^n$, then we have $\lambda_j=0$, which is a contradiction.
Therefore, in this case we have $B=\emptyset$ and $u\in \cA_0$, whence the 
statement $1)$ follows. 

If $P(0)\ne 0$, then we also have   
$B=\emptyset$ and $u\in \cA_0$, for $P_0(\lambda_j)=P(0)\ne 0$.  
Furthermore, since $P(0)\ne 0$, by Lemma \ref{Lma3.1}, $2)$ 
we have $0\not \in \mathcal Z_\Lambda(P)$ and $\cA_0 \cap \Ker P(D)=\{0\}$, 
whence $u=0$. Contradiction. Therefore $\rad(\Ker P(D))$ in this case contains no nonzero element 
and statement $2)$ follows. 
\epfv

Next, we show that Theorem \ref{Thm3.1} with some extra conditions holds also for commuting locally finite $R$-derivations. Recall that an $R$-derivation $D$ of an $R$-algebra $\cA$ is {\it locally finite (over $R$)} if  
for each $u\in \cA$, the $R$-submodule of $\cA$ 
spanned by elements $D^k u$ $(k\ge 0)$ over $R$ 
is finitely generated as an $R$-module.

\begin{propo}\label{NonCommLF-Case} 
Assume that $R$ is an integral domain of characteristic zero and 
$\cA$ is a reduced $R$-algebra that is torsion-free as an $R$-module. 
Denote by $K_R$ the field of fractions of $R$ and 
$\bar K_R$ the algebraic closure of $K_R$. 
Let $P(\xi)\in R[\xi]$ and  
$D=(D_1, D_2, \dots, D_n)$ be $n$ commuting 
locally finite $R$-derivations of $\cA$. 
Write $P(\xi)=\sum_{k=0}^d P_k(\xi)$ with $P_k(\xi)$ $(0\le k\le d)$ 
being homogeneous of degree $k$. 
Then the following statements hold:
\begin{enumerate}
\item[$1)$] if $P(0)= 0$ and  $P_k(\xi)$ $(1\le k\le d)$ have no nonzero common zeros in $\bar K_R^n$, then we have  
$ 
\rad(\Ker P(D))\subseteq \cA_0,
$ 
where $\cA_0=\bigcap_{i=1}^n 
(\sum_{m=1}^\infty \Ker D_i^m)$;
\item[$2)$] if $P(0)\ne 0$, then $\rad(\Ker P(D))=\{0\}$, 
and $\Ker P(D)$ is a MS of $\cA$.
\end{enumerate} 
\end{propo}

\pf Set $\bar\cA=\bar K_R\otimes_R\cA$.  
Since $\cA$ is torsion-free as an $R$-module, 
the standard map $\cA\simeq R\otimes_R\cA \to K_R\otimes_R \cA$ 
is injective, for by \cite[Prop. 3.3]{AM} 
$K_R\otimes_R \cA$ is isomorphic to the localization $S^{-1}\cA$ 
with $S=R\backslash\{0\}$. Since every field is absolutely flat, 
the standard map $K_R\otimes_R\cA \to \bar K_R\otimes_R \cA$ is also 
injective. Therefore, we may view $\cA$ as 
an $R$-subalgebra of $\bar\cA=\bar K_R\otimes_R\cA$ in the standard way and extend $D$ $\bar K_R$-linearly to $\bar\cA$, 
which we denote by $\bar D=(\bar D_1, \bar D_2, \dots, \bar D_n)$. 

Note that $\bar D_i$ $(1\le i\le n)$ are $n$ commuting  
$\bar K_R$-derivations of $\bar\cA$, 
which are also locally finite 
over $\bar K_R$. 
Then $\bar \cA$ by \cite[Proposition $1.3.8$]{E}) 
is decomposable w.r.t. $\bar D_i$ 
for each $1\le i\le n$. 
By applying Theorem \ref{Thm3.1} to $P(\bar D)$ and 
using the fact $\bar\cA_0 \cap \cA=\cA_0$ we see that  
the proposition follows.   
\epfv 
 
Next, we use the proposition above to show that 
Corollary \ref{Corol-2.4} with some extra conditions can be 
extended to the ordinary differential operators 
of some noncommutative algebras.

\begin{propo}\label{NonCommOneVariaCase} 
Let $R$, $\cA$ be as in Proposition \ref{NonCommLF-Case} and let  
$D$ be an arbitrary (single) $R$-derivation of $\cA$. 
Then for every univariate polynomial in 
$0\ne P(\xi)\in R[\xi]$,    
 the following statements hold: 
\begin{enumerate}
  \item[$1)$] if $P(0)=0$, then   we have 
$ 
\rad(\Ker P(D))\subseteq \rad(\cA_0),
$  
where $\cA_0=\sum_{j=1}^\infty \Ker D^j$; 
\item[$2)$] if $P(0)\ne 0$, then $\rad(\Ker P(D))=\{0\}$, 
and $\Ker P(D)$ is a MS of $\cA$.
\end{enumerate} 
\end{propo}

\pf The case $\deg P(\xi)=0$ is trivial. So we assume $\deg P(\xi)\ge 1$. 
Let $K_R$ be the field of fractions 
of $R$ with the algebraic closure $\bar K_R$, and set $\bar\cA=\bar K_R\otimes_R\cA$. 
As pointed out in the proof of Proposition \ref{NonCommLF-Case} 
we may view $\cA$ as an $R$-subalgebra of $\bar\cA$ in the standard way and 
extend $D$ $\bar K_R$-linearly to a $\bar K_R$-derivation of $\bar\cA$, 
which we denote by $\bar D$.  

Let $V=\Ker P(D)$. Then $V$ is an $R$-subspace 
of $\cA$ preserved by $D$.  
Set $\bar V=\bar K_R\otimes_R V$. Then $\bar D\,|_{\bar V}$ as a 
$\bar K_R$-linear map from $\bar V$ to $\bar V$ is algebraic over 
$\bar K_R$, for $P(D\,|_V)=P(D)\,|_V=0$ and hence $P(\bar D\,|_{\bar V})=0$. 
It is well-known (e.g.,  see \cite[Proposition $4.2$]{H}) 
that $\bar V$ can be decomposed as  a direct sum of the 
generalized eigenspaces of $\bar D\,|_{\bar V}$. 
Let $\bar \cB$ be the $\bar K_R$-subalgebra of $\bar \cA$ generated by 
elements of $\bar V$.  Then $\bar \cB$ is $\bar D$-invariant. Furthermore, 
by Eq.\,(\ref{LamMu-Rule}) it is easy to see that 
$\bar \cB$ is decomposable w.r.t. $\bar D\,|_{\bar \cB}$.


Now let $u\in \rad(\Ker P(D))$. Then there exists $N\ge 1$ such that  
$u^m\in \Ker P(D)$, and hence is also in $\bar \cB$, for all $m\ge N$.
Consequently, $u^m\in \rad(\Ker P(\bar D\,|_{\bar \cB}))$ for all $m\ge N$.
Write $P(\xi)=\sum_{k=0}^d P_k(\xi)$ (as before) with $P_k(\xi)$ $(0\le k\le d)$ 
homogeneous $k$ in $\xi$.
Then $P_k(\xi)$ $(1\le k\le d)$ have no nonzero common zero in $\bar K_R$, 
for $P(\xi)$ is a univariate polynomial of degree $d \ge 1$. 

If $P(0)=0$, then by applying Proposition \ref{NonCommLF-Case}, $1)$  
to $P(\bar D\,|_{\bar \cB})$ (as a differential operator of $\bar \cB$), 
we have $u^m\in \sum_{j=1}^\infty \Ker \bar D^j$ 
for all $m\ge N$. Since $\Ker D^j=\cA\cap\Ker \bar D^j$ for all $j\ge 1$, 
we further have $u^m\in \cA_0=\sum_{j=1}^\infty D^j$ 
for all $m\ge N$. Hence $u\in \rad(\cA_0)$ and statement $1)$ follows. 

If $P(0)\ne 0$, then by applying Proposition \ref{NonCommLF-Case}, $2)$ 
to $P(\bar D\,|_{\bar \cB})$ (as a differential operator of $\bar \cB$) 
we have $u^m=0$, and hence $u=0$, for $\cA$ is reduced.
Therefore statement $2)$ also holds. 
\epfv

We end this section with the following open problem which we believe is 
worthy of further investigations. 

\begin{oprob}
Let $R$ be an arbitrary commutative ring and $\cA$ an arbitrary unital noncommutative 
$R$-algebra. Let $D=(D_1, D_2, \dots, D_n)$ be $n$ $R$-derivations of $\cA$  
and let $Q(\xi)$ be a polynomial in $n$ noncommutative free variables 
$\xi=(\xi_1, \xi_2, \dots, \xi_n)$ over $R$. Set $a_0\!:=Q(0)$ and denote by $\ann_\ell(a_0)$  
the set of all elements $b\in \cA$ such that $a_0b=0$. Is it always true that $\rad(\Ker Q(D))\subseteq \rad\big(\ann_\ell(a_0)\big)$?
\end{oprob}

\section{\bf Some Applications to Locally Algebraic Derivations}\label{S4}
In this section we use some results proved in the last two 
sections to derive some properties of locally algebraic derivations 
and locally integral derivations.  

\begin{defi}\label{Def-LA}
Let $R$ be a unital commutative ring, $\cA$ 
an $R$-algebra and $D$ an $R$-derivation of $\cA$.
\begin{enumerate}
  \item[$1)$] We say $D$ is {\it algebraic over $R$} if  
there exists a nonzero polynomial $p(t)\in R[t]$ 
such that $p(D)=0$.

\item[$2)$] We say  $D$ is {\it locally algebraic over $R$} 
if for each $a\in \cA$, there exists a $D$-invariant $R$-subalgebra $\cA_1$ 
of $\cA$ containing $a$, and 
a nonzero polynomial $p_a (t)\in R[t]$ 
such that $p_a (D)\big|_{\cA_1}=0$. 
\end{enumerate}
\end{defi}

If $p(t)$ in statement $1)$ (resp., $p_a(t)$ in statement $2)$ for all $a\in \cA$) 
of the definition above can be chosen to be a monic polynomial, we say $D$ is 
{\it integral} (resp., {\it locally integral}) over $R$. 

An example of a derivation that is locally algebraic  
but not algebraic is as follows.

\begin{exam}
Let $x_i$ $(i\ge 1)$ be a sequence of free commutative variables and 
$\bC[x_i\,|\, i\ge 1]$ the polynomial algebra over $\bC$ in $x_i$ $(i\ge 1)$. Let 
$I$ be the ideal generated by $x_i^{i+1}$ $(i\ge 1)$ and $\cA=\bC[x_i\,|\, i\ge 1]/I$.
Then it can be readily verified that $D\!:=\sum_{i=1}^\infty x_i\p/\p x_i$ 
is a well-defined $\bC$-derivation of $\cA$, which is locally algebraic 
but not (globally) algebraic over $\bC$.     
\end{exam}

\begin{theo}\label{CommLocIntCase}
Let $R$ be a commutative ring and $\cA$ a commutative $R$-algebra 
such that the abelian group $(\cA, +)$ is torsion-free. 
Then for every $R$-derivation $D$ of $\cA$, which is locally integral over $R$, 
the image $\im D\!:=D(\cA)\subseteq \nil(\cA)$, where $\nil(\cA)$ denotes 
the nil-radical of $\cA$, i.e., the set of nilpotent elements of $\cA$.  
\end{theo}

\pf Let $a\in \cA$, and let $\cA_1$ be a $D$-invariant 
$R$-subalgebra of $\cA$ and $p_a(t)$ a monic polynomial in $R[t]$ 
such that $a\in \cA_1$ and $p_a(D)\,|_{\cA_1}=0$. 
Then by Corollary \ref{New-Corol-2.4} we have $d!(Da)^d=0$, where $d=\deg p_a(t)$. 
Since $(\cA, +)$ as an abelian group is torsion-free, we have $(Da)^d=0$, 
whence $Da\in \nil(\cA)$ and the theorem follows.  
\epfv

Since every nilpotent $R$-derivation of $\cA$ 
is locally integral over $R$, by Theorem \ref{CommLocIntCase} 
we immediately have the following:  

\begin{corol}\label{CommNilCase}
Let $R$, $\cA$ be as in Theorem \ref{CommLocIntCase} 
and let $D$ be a nilpotent $R$-derivation of $\cA$. 
Then $\im D\subseteq \nil(\cA)$.  
\end{corol}
 
Furthermore, from the proof of Theorem \ref{CommLocIntCase} it is also easy to see that we have the following:  

\begin{corol}\label{CommNilCase-2}
Let $R$ and $\cA$ be as in Theorem \ref{CommLocIntCase}. 
Assume further that $\cA$ is torsion-free   
as an $R$-module. Then for every $R$-derivation $D$ of $\cA$,  
which is locally algebraic over $R$, we have 
$\im D \subseteq \nil(\cA)$.  
\end{corol}

Next, we consider the $R$-derivations of some reduced $R$-algebra $\cA$ 
(not necessarily commutative), which are locally algebraic over $R$. 
%
%
%
 
\begin{theo}\label{NoncommDomainCase-1} 
Let $R$ be a unital integral domain of characteristic zero and let $\cA$ be 
a unital reduced $R$-algebra (not necessarily commutative) 
that is torsion-free as an $R$-module. 
Then $\cA$ has no nonzero $R$-derivations that 
are locally algebraic over $R$. In particular, 
$\cA$ has no nonzero nilpotent $R$-derivations. 
\end{theo}

\pf Let $D$ be an $R$-derivation of $\cA$ that is locally algebraic over $R$. 
Let $a\in \cA$, and $\cA_1$ be a $D$-invariant 
$R$-subalgebra of $\cA$ and $0\ne p_a(t)\in R[t]$   
such that $a\in \cA_1$ and $p_a(D)\,|_{\cA_1}=0$. 
Then $a^m\in \cA_1\subseteq \Ker \, p_a(D)$ for all $m\ge 1$, 
whence $a\in\rad(\Ker \, p_a(D))$. 

Replacing $p_a(t)$ by $tp_a(t)$ we assume $p_a(0)=0$. 
Then by applying Proposition \ref{NonCommOneVariaCase}, $1)$
to the differential operator $p_a(D)$, we have 
$a\in \rad(\cA_0)$, where $\cA_0=\sum_{i=1}^\infty \Ker D^i$.
Consequently, $\rad(\cA_0)=\cA$. Then by \cite[Lemma 2.4]{MS} we have 
$\cA_0=\cA$, i.e., $D$ is locally nilpotent. 

Let $K_R$ be the field of fractions of $R$ and $\cB\!:= K_R\otimes_R\cA$. 
As pointed out in the proof of Proposition \ref{NonCommLF-Case},  
we may view $\cA$ as an $R$-subalgebra 
of $\cB$ in the standard way and 
extend $D$ $K_R$-linearly to a $K_R$-derivation of $\cB$, 
which we denote by $\bar D$.  

Let $a$, $p_a(t)$ be fixed as above, and $N\ge 1$ such that $D^Na=0$. 
Write $p_a(t)=t^kh(t)$ for some $k\ge 1$ and $h(t)\in K_R[t]$ 
with $h(0)\ne 0$. Then $p_a(\bar D)a=0$ and $\bar D^Na=0$. Since  
$\gcd\,(p_a(t),\, t^N)=t^\ell$ in $K_R[t]$ with $\ell=\min\{k, N\}$, 
we have $\bar D^\ell a=0$. Hence $D^ka=\bar D^k a=0$. 
Since $a$ is an arbitrary element of $\cA$, we have $D^k=0$.   
Then by \cite[Lemma $6.1$]{Open-LFNED} we have $D=0$, whence 
the theorem follows. 
\epfv

One remark on Theorem \ref{NoncommDomainCase-1} is that, without 
the characteristic zero condition, the theorem   
may be false, which can be seen from the following example. For more integral derivations 
of algebras over a field of characteristic $p>0$, see \cite{N}.

\begin{exam}
Let $K$ be a field of characteristic $p>0$, $\cA=K[x]$ 
and $D=d/dx$. Then $D^p=0$. Hence $D$ is a nonzero 
$K$-derivation of $\cA$ that is algebraic over $K$. 
\end{exam}

One immediate consequence of Theorem \ref{CommLocIntCase}, Corollary \ref{CommNilCase-2} and 
Theorem \ref{NoncommDomainCase-1} is the following corollary which in some sense gives an affirmative answer 
to the so-called LNED conjecture proposed in \cite{Open-LFNED}
for nilpotent, or locally integral, or locally algebraic derivations 
of certain algebras.

\begin{corol} 
$1)$ Let $R$, $\cA$ be as in Theorem \ref{CommLocIntCase} and let $D$ be an $R$-derivation of $\cA$. If $D$ is  
locally integral over $R$, then $D$ maps every $R$-subspace of $\cA$ to a MS of $\cA$.

$2)$ Let $R$, $\cA$ be as in Corollary \ref{CommNilCase-2} 
and let $D$ be an $R$-derivation of $\cA$. If $D$ is  
locally algebraic over $R$, then $D$ maps every $R$-subspace of $\cA$ 
to a MS of $\cA$.
\end{corol}

We end this section with the following proposition which is not needed elsewhere 
in this paper but is interesting on its own for the study of 
the radical of the kernel of a derivation. 

\begin{propo}\label{Propo-4.9}
Let $R$ be a commutative ring and $\cA$ a 
reduced $R$-algebra (not necessarily commutative) 
such that $(\cA, +)$ is torsion-free.  
Let $r\ge 1$, $a\in \cA$ and $D$ be an $R$-derivation of $\cA$ 
such that $D^r a^m=0$ for all $1\le m\le 2^{r-1}$.  
Then $a\in \Ker D$. Consequently, $\Ker P(D)\subseteq \rad(\Ker D)= \rad(\Ker D^k)$ for all $k\ge 1$. 
\end{propo}

Note that when $\cA$ is commutative, the lemma follows easily from 
Theorem \ref{CommCase} with $P(D)=D^r$. 
Here we give a proof without the commutativity of $\cA$. \\

\underline{\bf Proof of Proposition \ref{Propo-4.9}}: We use induction on $r$. 
The case $r=1$ is obvious. 
So assume $r\ge 2$. Then $2r-2\ge r$ and  
for each $1\le k\le 2^{r-2}$, by the Leibniz rule we have 
$$
0=D^{2r-2} a^{2k}=\sum_{i=0}^{2r-2}\binom{2r-2}i 
(D^i a^k)(D^{2r-2-i}a^k).
$$
Since $D^j a^k=0$ for all $j\ge r$, there is only one term in 
the sum above that may not be equal to $0$, namely, 
the term with $i=r-1$. Therefore  
$\binom{2r-2}{r-1} (D^{r-1} a^k)^2=0$. 
Since $\cA$ is reduced and $(\cA, +)$ 
is torsion-free, we have $D^{r-1} u^k=0$ for all $1\le k\le 2^{(r-1)-1}$. 
Then by the induction assumption we have $a\in\Ker D$.

Now let $k\ge 1$ and $u\in \rad(\Ker D^k)$. Then there exists $N\ge 1$ 
such that $u^m\in \Ker D^k$ for all $m\ge N$. Applying the result proved  
above to $u^m$ $(m\ge N)$ we have $u^m\in \Ker D$ for all $m\ge N$. 
Hence $u\in \rad(\Ker D)$, and $\rad(\Ker D^k)\subseteq \rad(\Ker D)$. 
Conversely, since $\Ker D$ as 
an $R$-subalgebra of $\cA$ is closed under the multiplication  
and $\Ker D \subseteq \Ker D^k$, we also have 
$\Ker D \subseteq \rad(\Ker D)$ and 
$\rad(\Ker D) \subseteq \rad(\Ker D^k)$.  
Hence the proposition follows. 
\epfv

\renewcommand{\theequation}{\thesection.\arabic{equation}}
\renewcommand{\therema}{\thesection.\arabic{rema}}
\setcounter{equation}{0}
\setcounter{rema}{0}

\section{\bf A Differential Vandemonde Determinant}\label{S5}
Throughout this section {\it $\cA$ stands for 
a commutative ring and $D$ for a derivation of $\cA$}.

\begin{propo}\label{DiffVandeMonde}
Let $\cA$ and $D$ be fixed as above. 
Then for all $f\in \cA$ and $n\ge 1$, we have
\begin{align}\label{DiffVandeMonde-eq1}
\det\,\begin{pmatrix}
f & f^2 & \cdots & f^n \\
D(f) & D(f^2) & \cdots & D(f^n) \\
D^2(f) & D^2(f^2) & \cdots & D^2(f^n) \\
 \vdots&\vdots&{}&\vdots \\
D^{n-1}(f) & D^{n-1}(f^2) & \cdots & D^{n-1}(f^n) 
\end{pmatrix}
=\alpha_n (Df)^{\frac12 n(n-1)} f^n,
\end{align}
where $\alpha_n=\prod_{k=1}^{n-1} k!$.
\end{propo}

The main idea of the proof of the proposition above is 
to show that the matrix in Eq.\,(\ref{DiffVandeMonde-eq1}) 
can be transformed by some elementary column 
operations to an upper triangular matrix whose  
 $i$-{th} diagonal entry is equal to 
$(i-1)!(Df)^{i-1}f$ for all $1\le i\le n$. 
For example, for the case $n=2$, by subtracting from the second column
the multiple of the first column by $f$ we get  
\begin{align} \label{DiffVandeMonde-peq0}
\begin{pmatrix}
f&f^2\\
D(f)&D(f^2)
\end{pmatrix}
\Longrightarrow
\begin{pmatrix}
f&0 \\
D(f)&fD(f)
\end{pmatrix}.
\end{align}

To see that this can be achieved for all $n\ge 2$, it suffices to show the following lemma, 
from which Proposition \ref{DiffVandeMonde} immediately follows. 

\begin{lemma}\label{Lma-5.2} 
Let $D$ and $f$ be as in Proposition \ref{DiffVandeMonde} and $k\ge 2$. Then 
there exist $\alpha_{k,j}\in \cA$ $(1\le j\le k-1)$ such that for each $0\le i\le k-1$, we have
\begin{align}\label{Lma-5.2-eq1}
D^i(f^k)-\sum_{j=1}^{k-1}\alpha_{k,j} f^{k-j} D^i(f^j)=\delta_{i, k-1} (k-1)!(Df)^{k-1}f,
\end{align}
where $\delta_{i, k-1}$ is the Kronecker delta function. 
%
\end{lemma}
\pf We use induction on $k$. If $k=2$, then $\alpha_{2,1}=1$ 
 solves the equations in Eq.\,(\ref{Lma-5.2-eq1}), as already pointed out in Eq.\,(\ref{DiffVandeMonde-peq0}) above. 

Assume that the lemma holds for some $k\ge 2$ and consider the case $k+1$. 
By writing $f^{k+1}$ as $f\cdot f^k$ and applying 
the Leibniz rule, we have for each $0 \le i\le k$  
\begin{align*}
D^i&(f^{k+1})-fD^i(f^k)=\sum_{\ell =0}^{i-1} \binom i\ell 
(D^{i-\ell} f) D^\ell(f^k)
\intertext{Applying the induction assumption to $D^\ell(f^k)$ and 
noticing that $\ell=k-1$, if and only if $i=k=\ell+1$, since $\ell\le i-1\le k-1$:}
&=\delta_{i, k} k(k-1)!(Df)^{k}f  +
\sum_{\ell =0}^{i-1} \binom i\ell 
(D^{i-\ell} f)
\left( \sum_{j=1}^{k-1}\alpha_{k,j} f^{k-j} D^\ell(f^j) 
\right)\\
&=\delta_{i, k} k!(Df)^{k}f  +
\sum_{j=1}^{k-1}\alpha_{k,j} f^{k-j}
\sum_{\ell =0}^{i-1} \binom i\ell 
(D^{i-\ell} f) D^\ell(f^j) \\
&=\delta_{i, k} k!(Df)^{k}f  + \sum_{j=1}^{k-1}\alpha_{k,j} f^{k-j}
 \Big(D^i(f^{j+1}) -fD^i(f^j)\Big) \\
&=\delta_{i, k} k!(Df)^{k}f  +\alpha_{k,k-1} fD^i(f^k)-\alpha_{k,1} f^{k}D(f)\\
&\hspace{45mm} +
\sum_{j=2}^{k-1}(\alpha_{k, j-1}-\alpha_{k, j}) f^{k+1-j} D^i(f^j) . 
\end{align*}

Set 
\begin{align}\label{RecurRelation}
\alpha_{k+1,j}\!:=
\begin{cases}
-\alpha_{k,1} &\text{ if  } \,\, j=1;\\
\alpha_{k, j-1}-\alpha_{k,j} &\text{ if  } \,\,  2\le j\le k-1;\\
1+\alpha_{k, k-1} &\text{ if  } \,\,  j=k. 
\end{cases}
\end{align}
Then $\alpha_{k+1,j}$ $(1\le j\le k)$ solve the equations in 
Eq.\,(\ref{Lma-5.2-eq1}) for the case $k+1$ and $0\le i\le k$. 
Hence, by induction the lemma follows. 
\epfv

\begin{rmk}\label{OtherApp}
One application of the formula in Eq.\,(\ref{DiffVandeMonde-eq1}) is as follows.  
We first apply the formula to some special function $f(x)$ and derivation $D$, 
and then evaluate $x$ at a fixed point $c$. By doing so, we may 
get formulas for the determinants of several families of matrices, e.g.,  
letting $f=x^d$, $D=x^m \frac{d}{dx}$ and $c=\pm 1$ for all $d, m \in \bZ$.
In particular, if we choose $d=-1$, $m=0$ and $c= 1$, then with a little more argument  
we get the following formula with $0!\!:=1$ for all $n\ge 1$: 
\begin{align}
\det\big( (i+j-2)!\big)_{1\le i, j\le n}&=\big(\prod_{k=1}^{n-1} k! \big)^2. 
\end{align}
%
\end{rmk}

Another consequence of Proposition \ref{DiffVandeMonde} is 
the following:  

\begin{propo}\label{Propo-5.4}
Let $D$ be a derivation of $\cA$, $\xi$ be a free variable and  
$0\ne P(\xi)=\sum_{i=0}^d c_i\xi^i\in \cA[\xi]$. Let $f\in\cA$ be   
such that $f^m\in\Ker P(D)$ for all $1\le m\le d+1$. Then   
for each $0\le i\le d$,  we have  
\begin{align}
 \alpha_{d+1} c_i (Df)^{\frac12 d(d+1)} f^{d+1}=0,
\end{align} 
where $\alpha_{d+1}=\prod_{k=1}^{d} k!$
\end{propo}

\pf Let $B$ be the transpose of the matrix in Eq.\,(\ref{DiffVandeMonde-eq1}) 
with $n-1=d$. Since $P(D)(f^m) =0$ for all $1\le m\le d+1$, we have $Bv=0$, 
where $v$ denotes the column vector $(c_0, c_1, \dots, c_d)^\tau$.  
Then  the proposition follows from Eq.\,(\ref{DiffVandeMonde-eq1}). 
\epfv

\begin{corol}
Let $D$, $f$, $P(\xi)$ be as in Proposition \ref{Propo-5.4}. 
Assume further that $(\cA, +)$ is torsion-free and that 
$c_i$ is not zero nor a zero-divisor of $\cA$ for some $0\le i\le d$.  
Then $fDf$ is nilpotent.
\end{corol}

\pf By Proposition \ref{Propo-5.4} we have 
$\alpha_{d+1} c_i (Df)^{\frac12 d(d+1)} f^{d+1}=0$, and hence, 
$f^{d+1} (Df)^{\frac12 d(d+1)} =0$, for 
$(\cA, +)$ is torsion-free and $c_i$ is not zero nor 
a zero-divisor of $\cA$.   
Then $(fDf)^m=0$ for all $m\ge \max\{d+1,\,\frac12 d(d+1)\}$, whence  
the corollary follows.
\epfv


\begin{thebibliography}{FLM2}

 \bibitem[AM]{AM} M. F. Atiyah and I. G. Macdonald, {\it Introduction to Commutative Algebra.}
 Addison-Wesley Publishing Co., 1969. [MR0242802].

\bibitem[BCW]{BCW} H. Bass, E. Connell and D. Wright, {\it The Jacobian Conjecture, Reduction of Degree and Formal Expansion of the Inverse}. Bull.  Amer. Math.  Soc.  \textbf{7}, (1982), 287--330. 


\bibitem[DEZ]{DEZ} H. Derksen, A. van den Essen and W. Zhao, {\it The Gaussian Moments Conjecture and the Jacobian Conjecture}. To appear in {\it Israel J. Math.}. See also arXiv:1506.05192 [math.AC]. 


\bibitem[E1]{E} A. van den Essen, {\it Polynomial Automorphisms and the Jacobian Conjecture}.
Prog. Math., Vol.190, Birkh\"auser Verlag, Basel, 2000. 

\bibitem[E2]{E2} A. van den Essen, 
{\it Introduction to Mathieu Subspaces}. 
``International Short-School/Conference on Affine Algebraic Geometry and the Jacobian Conjecture" at Chern Institute of Mathematics, Nankai University, Tianjin, China. July 14-25, 2014.


\bibitem[EH]{EH} A. van den Essen and L. C. van Hove, {\it Mathieu-Zhao spaces of polynomial rings}. arXiv:1907.06106 [math.AC].

\bibitem[EKC]{EKC} A. van den Essen, S. Kuroda and A. J. Crachiola, {\it Polynomial Automorphisms 
and the Jacobian Conjecture: New Results from the Beginning of the 21st
Century}. Frontiers in Mathematics. Birkh\"auser $2021$. 

\bibitem[EN]{EN} A. van den Essen and S. Nieman, {\it Mathieu-Zhao Spaces of Univariate Polynomial Rings with Non-zero Strong Radical}.  J. Pure Appl. Algebra, {\bf 220} (2016), no.\,9, 3300--3306.  

\bibitem[EZ1]{EZ1} A. van den Essen and W. Zhao, {\it Mathieu Subspaces of Univariate Polynomial Algebras}. J. Pure Appl. Algebra. 217 (2013), no.7, 1316–-1324. See also arXiv:1012.2017[math.AC].

\bibitem[EZ2]{EZ2} A. van den Essen and W. Zhao, {\it On the Image Conjecture for Locally Finite Derivations and $\mathcal E$-Derivations}. 
 J. Pure Appl. Algebra 223 (2019), no.\,{\bf 4}, 1689-1698. See also arXiv:1708.05813 [math.AC].  

\bibitem[EWZ]{EWZ} A. van den Essen, D. Wright and W. Zhao, 
{\it Images of Locally Finite Derivations of Polynomial Algebras in Two Variables}. J. Pure Appl. Algebra {\bf 215} (2011), no.9, 2130-2134. [MR2786603]. See also arXiv:1004.0521[math.AC].

\bibitem[Hu]{H} J. E. Humphreys, (1972), 
{\it Introduction to Lie Algebras and Representation Theory}.  
Graduate Texts in Mathematics, Springer, 1972. 

\bibitem[Ke]{K} O. H. Keller, 
{\it Ganze Gremona-Transformationen}. Monats. Math. Physik {\bf 47} (1939), no.\,1, 299-306. 

\bibitem[Kh1]{Kharchenko1} V.K. Kharchenko, {\it Differential identities of prime rings}.  Algebra Logika 17 (1978) 220–238 (English translation in Algebra Logic 17 (1978) 154–168). [MR0541758] (81f:16025).

\bibitem[Kh2]{Kharchenko2} V.K. Kharchenko, {\it Automorphisms and Derivations of Associative Rings}. Math. Appl. (Soviet Ser.), vol. 69, Kluwer Academic Publishers Group, Dordrecht, 1991, translated from the Russian by L. Yuzina. [MR1174740] (93i:16048).

\bibitem[L]{L} S. Lefschetz, {\it Differential Equations: Geometric Theory}. Reprinting of the
second edition. Dover Publications, Inc., New York, 1977. [MR0435481].

\bibitem[M]{Ma} O. Mathieu, {\it Some Conjectures about Invariant Theory and Their Applications.} Alg\`ebre non commutative, groupes quantiques et invariants (Reims, 1995), 263--279, S\'emin. Congr., 2, Soc. Math. France, Paris, 1997. [MR1601155].

\bibitem[N]{N} A. Nowicki, {\it Integral Derivations}. J. Alg. {\bf 110} (1987), no.\,1, 262-276. 


%
 
 

\bibitem[Z1]{IC} W. Zhao, {\it Images of Commuting  Differential Operators of Order One with Constant Leading Coefficients}.  
J. Alg. {\bf 324} (2010),  no. 2, 231--247. [MR2651354].  
See also arXiv:0902.0210 [math.CV]. 


\bibitem[Z2]{GIC} W. Zhao, {\it Generalizations of the Image Conjecture and the Mathieu Conjecture}.  J. Pure Appl. Alg. {\bf 214} (2010), 1200-1216. See also arXiv:0902.0212 [math.CV].  
 

\bibitem[Z3]{MS} W. Zhao, {\it Mathieu Subspaces of Associative Algebras}. J. Alg.  
{\bf 350} (2012), no.2, 245-272. 
See also arXiv:1005.4260 [math.RA].
 
\bibitem[Z4]{Open-LFNED} W. Zhao, {\it Some Open Problems on Locally Finite or Locally Nilpotent Derivations and $\cE$-Derivations}. Commun. Contemp. Math. 20 (2018), no.\,{\bf 4}, 1750056, 25pp. See also arXiv:1701.05992 [math.RA]. 

\bibitem[Z5]{Idem} W. Zhao, {\it Idempotents in Intersection of the Kernel and the Image of Locally Finite Derivations and $\cE$-derivations}.  Eur. J. Math. 4 (2018), no.\,{\bf 4}, 1491-1504. See also arXiv:1701.05993 [math.RA]. 

\bibitem[Z6]{Alg-LFNED} W. Zhao, {\it The LFED and LNED Conjectures 
for Algebraic Algebras}. Linear Algebra Appl.\,{\bf 534} (2017), 181-194. 
  See also arXiv:1701.05990 [math.RA]. 

\bibitem[Z7]{LaurentPolyCase} W. Zhao, {\it The LFED and LNED Conjectures 
for Laurent Polynomial Algebras}. Under submission. See also arXiv:1701.05997 [math.AC].

\bibitem[Z8]{OneVariableCase} W. Zhao, {\it Images of Ideals under Derivations and $\cE$-Derivations of Univariate Polynomial Algebras  over a Field  of Characteristic Zero}. Under submission. See also arXiv:1701.06125 [math.AC]. 

\end{thebibliography}
\end{document}